\documentclass{article}
\usepackage{pdfpages}
\usepackage[framemethod=TikZ]{mdframed}
\usepackage{multicol}
\usepackage{tikz}
\usepackage{tikz-3dplot}
\usepackage{subcaption}
\usepackage{mathrsfs}
\usepackage{hyperref}
\usepackage{algpascal}
\usepackage{algorithmicx}
\usepackage{algorithm}
\usepackage{algpseudocode}
\usepackage{blindtext}
\usepackage{amssymb}
\usepackage{authblk}
\usepackage{natbib}
\usepackage[a4paper, total={5.5in, 9.5in}]{geometry}
\usepackage[utf8]{inputenc}
\usepackage[english]{babel}
\usepackage{amsmath}
\newtheorem{theorem}{Theorem}
\newtheorem{lemma}{Lemma}

\title{Finding Euclidean Distance to a Convex Cone Generated by a Large Number of Discrete Points}
\author[1]{Ali Fattahi}
\author[2]{Sriram Dasu}
\author[1]{Reza Ahmadi}
\affil[1]{Anderson School of Management, University of California, Los Angeles, Los Angeles, California 90095}
\affil[2]{Marshall School of Business, University of Southern California, Los Angeles, Los Angeles, California 90089}
\date{}                     
\begin{document}
\maketitle
\begin{abstract}
In this paper, we study the problem of finding the Euclidean distance to a convex cone generated by a set of discrete points in $\mathbb{R}^n_+$. In particular, we are interested in problems where the discrete points are the set of feasible solutions of some binary linear programming constraints. This problem has applications in manufacturing, machine learning, clustering, pattern recognition, and statistics. Our problem is a high-dimensional constrained optimization problem. We propose a Frank-Wolfe based algorithm to solve this non-convex optimization problem with a convex-noncompact feasible set. Our approach consists of two major steps: presenting an equivalent convex optimization problem with a non-compact domain, and finding a compact-convex set that includes the iterates of the algorithm. We discuss the convergence property of the proposed approach. Our numerical work shows the effectiveness of this approach.

\textbf{Subject classifications:} Frank-Wolfe; Configurable Products; Clustering; Machine Learning;
\end{abstract}

\section{Introduction}

Given a large set of discrete points in $\mathbb{R}^n_+$, denoted by $\mathbb{Y}$, in this paper, we are interested in obtaining the Euclidean distance from a target point to the convex cone generated by $\mathbb{Y}$, where $n$ is in the order of thousands and $m$ is significantly large. In particular, our approach is suitable for problems where $\mathbb{Y}$ is the set of feasible solutions of some Binary Linear Programming (BLP) constraints, for which we have $m=\mathcal{O}(2^n)$. This problem is represented as minimizing a non-convex function over a convex-noncompact domain.

Our problem is motivated by an application in a large auto manufacturer, in which, $n$ represents options that define the car configurations, the target point is an estimation of future demands of options, and the convex cone constitutes the set of all feasible points---i.e., points that show a producible set of car configurations \citep{FattahiEtAll2016,FattahiEtAll2017}.  This problem can generally be found in all manufacturing systems where products are configured based on a set of options available for the customers. Other applications include variants of classical clustering problem. In some clustering problems, the objective is to find the Euclidean distance of a point to the convex hall of a set of points in $\mathbb{R}^n$, called a cluster. Whereas, our problem is to find the Euclidean distance of a target point to the convex cone of a set of points in $\mathbb{R}^n_+$. If $\mathbb{Y}$ is given numerically, our problem can be formulated as a non-negative least squares problem (see, for example, \citet{FrancEtAl2005,BoutsidisAndDrineas2009,Potluru2012}).

The well-known Frank-Wolfe Algorithm (FWA) \citep{FrankAndWolfe1956}, a.k.a. conditional gradient algorithm, has been recently used to solve high-dimensional constraint optimization problems in machine learning, pattern recognition, clustering, and statistics. The original FWA minimizes a convex function over a convex-compact domain. \citet{Jaggi2013,Lacoste_JulienAndJaggi2015} have looked at the algorithmic variants of the FWA, and  \citet{ReddiEtAl2016,LafondEtAl2016,HazanAndKale2012} have recently extended the application of the FWA to minimizing a non-convex objective function and stochastic optimization.

In this paper, we apply the FWA to a nonconvex optimization problem with a convex-noncompact domain. We first show that our problem is equivalent to a problem with a convex objective function and a convex-noncompact feasible set. Then, we show that there exists a compact-convex set that contains the iterates of the FWA, and we characterize the diameter of this compact set. Consequently, our proposed Frank-Wolfe based approach solves this problem at a linear convergence rate, i.e. the optimization error after $k$ iterations will decrease with $\mathcal{O}(\frac{1}{k})$.

\section{Problem Description} \label{SectionProblemDescription}

Let $\mathbb{Y}=\{\mathbf{y}^1,\mathbf{y}^2,\dots,\mathbf{y}^m\}$, where $\mathbf{y}^i\in\mathbb{R}^n_+$, for all $i$. We assume in this paper that there exists $\mathbf{y}\in\mathbb{Y}$ such that $\mathbf{y}\ne\mathbf{0}$. Define $\mathbb{X}:=conv(\mathbb{Y})$, and $\mathbb{Z}:=cone(\mathbb{X})$. A reference point $\hat{\mathbf{z}}\in \mathbb{R}^n_+$, $\hat{\mathbf{z}}\ne\mathbf{0}$, is given and we have to find a feasible point $\mathbf{z}^*\in\mathbb{Z}$ which has the closest Euclidean distance to the reference point $\hat{\mathbf{z}}$.
Our problem then is defined as
$\mathscr{P}': \min_{\mathbf{z}\in\mathbb{Z}}  \|\mathbf{z}-\hat{\mathbf{z}}\|_2^2$. This problem can be equivalently represented as $\mathscr{P}'': \min_{\mathbf{x}\in\mathbb{X},\lambda\ge0}  \|\lambda\mathbf{x}-\hat{\mathbf{z}}\|_2^2$.

\begin{figure}[h]
\begin{center}
    \begin{subfigure}[b]{.49\textwidth}
        \centering
    \tdplotsetmaincoords{70}{120}
    \begin{tikzpicture}[tdplot_main_coords, scale = 1.1]

    \draw[opacity=0.2,thin,dashed](0,0,0)--(0,4,0)--(0,4,4)--(0,0,4)--(0,0,0)--cycle;
    \draw[opacity=0.2,thin,dashed](0,0,0)--(4,0,0)--(4,4,0)--(0,4,0)--(0,0,0)--cycle;
    \draw[opacity=0.2,thin,dashed](0,0,0)--(4,0,0)--(4,0,4)--(0,0,4)--(0,0,0)--cycle;
    \draw[opacity=0.2,thin,dashed](4,0,0)--(4,4,0)--(4,4,4)--(4,0,4)--(4,0,0)--cycle;
    \draw[opacity=0.2,thin,dashed](0,4,0)--(4,4,0)--(4,4,4)--(0,4,4)--(0,4,0)--cycle;
    \draw[opacity=0.2,thin,dashed](0,0,4)--(4,0,4)--(4,4,4)--(0,4,4)--(0,0,4)--cycle;

    \draw[fill=black] (0,4,0) circle (0pt) node[right]{$\mathbb{R}^3_+$};
        \draw[thick,->](0,0,0)--(3,0,0) node[below left]{};
        \draw[thick,->](0,0,0)--(0,3,0) node[below right]{};
        \draw[thick,->](0,0,0)--(0,0,3) node[left]{};

        \draw[thick,->,blue](0,0,0)--(1,1,2) node[right]{$\mathbf{y}^1$};

        \draw[thick,->,blue](0,0,0)--(0,2,3) node[below right]{$\mathbf{y}^2$};

        \draw[thick,->,blue](0,0,0)--(2,1,3) node[left]{$\mathbf{y}^3$};

        \draw[thick,->,blue](0,0,0)--(3,0,2) node[below left]{$\mathbf{y}^4$};

        \draw[thick,->,blue](0,0,0)--(0,0,2) node[right]{$\mathbf{y}^5$};

        \draw[thick,->,brown](0,0,0)--(1,1,0) node[below right]{$\hat{\mathbf{z}}$};

    \end{tikzpicture}
    \end{subfigure}
    \begin{subfigure}[b]{.49\textwidth}
        \centering
    \tdplotsetmaincoords{70}{120}
    \begin{tikzpicture}[tdplot_main_coords, scale = 1.1]

    \draw[opacity=0.2,thin,dashed](0,0,0)--(0,4,0)--(0,4,4)--(0,0,4)--(0,0,0)--cycle;
    \draw[opacity=0.2,thin,dashed](0,0,0)--(4,0,0)--(4,4,0)--(0,4,0)--(0,0,0)--cycle;
    \draw[opacity=0.2,thin,dashed](0,0,0)--(4,0,0)--(4,0,4)--(0,0,4)--(0,0,0)--cycle;
    \draw[opacity=0.2,thin,dashed](4,0,0)--(4,4,0)--(4,4,4)--(4,0,4)--(4,0,0)--cycle;
    \draw[opacity=0.2,thin,dashed](0,4,0)--(4,4,0)--(4,4,4)--(0,4,4)--(0,4,0)--cycle;
    \draw[opacity=0.2,thin,dashed](0,0,4)--(4,0,4)--(4,4,4)--(0,4,4)--(0,0,4)--cycle;

    \draw[fill=black] (0,4,0) circle (0pt) node[right]{$\mathbb{R}^3_+$};
        \draw[thick,->](0,0,0)--(3,0,0) node[below left]{};
        \draw[thick,->](0,0,0)--(0,3,0) node[below right]{};
        \draw[thick,->](0,0,0)--(0,0,3) node[left]{};

        \filldraw[draw=brown,fill=brown,opacity=0.2](0,0,0)--(4,0,8/3)--(4,0,4)--(0,0,4)--(0,0,0)--cycle;
        \filldraw[draw=brown,fill=brown,opacity=0.23](0,0,0)--(0,8/3,4)--(0,0,4)--(0,0,0)--cycle;
        \filldraw[draw=brown,fill=brown,opacity=0.26](0,0,4)--(0,8/3,4)--(2,2,4)--(4,1,4)--(4,0,4)--(0,0,4)--cycle;
        \filldraw[draw=brown,fill=brown,opacity=0.29](0,0,0)--(0,8/3,4)--(2,2,4)--(0,0,0)--cycle;
        \filldraw[draw=brown,fill=brown,opacity=0.32](0,0,0)--(4,0,8/3)--(4,1,4)--(2,2,4)--(0,0,0)--cycle;
        \filldraw[draw=brown,fill=brown,opacity=0.35](4,0,8/3)--(4,1,4)--(4,0,4)--(4,0,8/3)--cycle;
        \draw[thick,->,brown](0,0,0)--(1,1,0) node[below right]{$\hat{\mathbf{z}}$};

        \draw[fill=black] (0,1,1.5) circle (0pt) node[right]{$\mathbb{Z}$};

    \end{tikzpicture}
    \end{subfigure}
    \hfill
    \caption{What is the Euclidean distance between $\hat{\mathbf{z}}$ and the convex cone generated by $\{\mathbf{y}^1,\mathbf{y}^2,\mathbf{y}^3,\mathbf{y}^4,\mathbf{y}^5\}$? (Note: $\mathbf{y}^1=(1,1,2),\ \mathbf{y}^2=(0,2,3),\ \mathbf{y}^3=(2,1,3),\ \mathbf{y}^4=(3,0,2),\ \mathbf{y}^5=(0,0,2)$, and $\hat{\mathbf{z}}=(1,1,0)$.)\label{RCPillustFig}}
\end{center}
\end{figure}
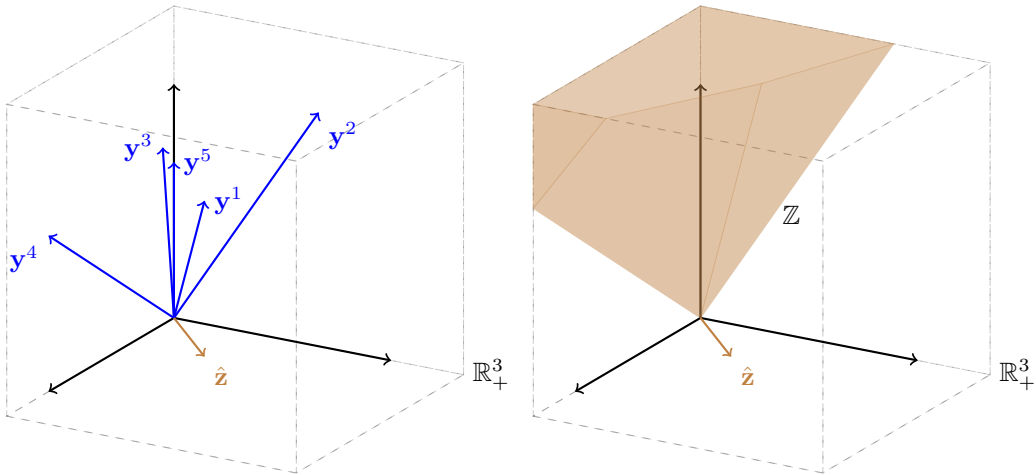

Fig. \ref{RCPillustFig} provides an illustration in $\mathbb{R}^3_+$. Given vectors $\mathbf{y}^1,\dots,\ \mathbf{y}^5$, we want to find the distance between a target point, $\hat{\mathbf{z}}$, and the convex cone generated by $\mathbf{y}^1,\dots,\ \mathbf{y}^5$.

We are particularly interested in the instances where $n$ is in the order of thousands and $m=\mathcal{O}(2^n)$, when  $\mathbb{Y}$ is the set of feasible solutions to a BLP, in which case $n$ will be the number of binary variables. This problem is very difficult. Problem $\mathscr{P}'$ has a convex objective function, while $\mathscr{P}''$ has a non-convex objective function. Moreover, problem $\mathscr{P}'$ involves nonlinearity in $\mathbf{z}\in\mathbb{Z}$ because of the definition of cone.

\begin{figure}[h]
  \centering
  \includegraphics[width=12cm, trim= 0cm 2cm 0cm 0.7cm]{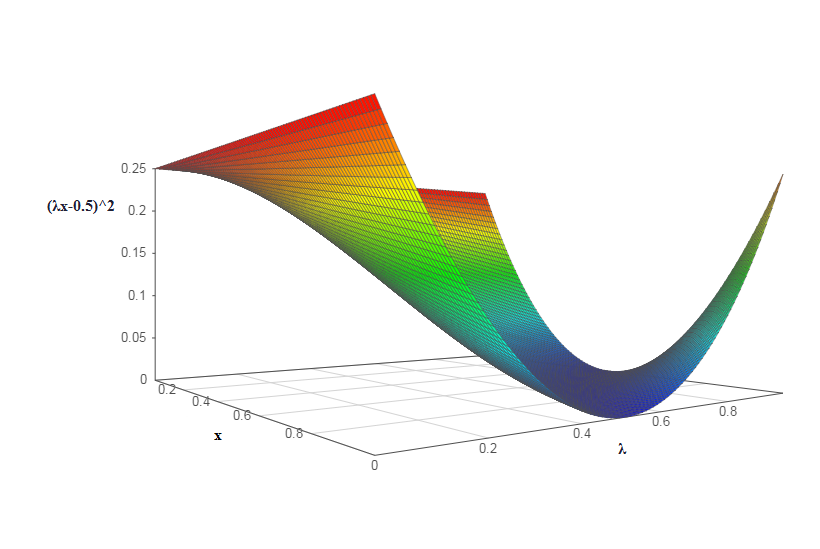}
  \caption{The graph of $(\lambda x-0.5)^2$ over the domain $0\le x\le1,\ 0\le\lambda\le1$}\label{PnotConvFig}
\end{figure}

\begin{lemma}\label{PnotConvexProblemLemma}
$\|\lambda \mathbf{x}-\hat{\mathbf{z}}\|_2^2$ is not necessarily convex over the domain $\mathbf{x}\in \mathbb{X},\ \lambda\ge0$.
\end{lemma}
\noindent\textsc{Proof.} Consider an instance of $\mathscr{P}''$ with $n=1$, $\mathbb{Y}=\{0,1\}$, and $\hat{z}=0.5$. For this instance, $\mathscr{P}''$ can be written as $\min_{0\le x\le1,\ \lambda\ge0}(\lambda x-0.5)^2$. Fig. \ref{PnotConvFig} shows the graph of $(\lambda x-0.5)^2$ over the domain $0\le x\le1,\ 0\le\lambda\le1$. It is seen that the objective function is not convex (consider, for example, the diagonal from (0,0) to (1,1)).
$\Box$

A possible approach is to formulate the feasible region of $\mathscr{P}'$ as a set of mixed-integer nonlinear programming (MINLP) constraints. The optimal solution $\mathbf{z}^*$ may be in the interior or boundary of $\mathbb{Z}$. We remark that $\mathbf{z}^*$ can be represented as a non-negative combination of at most $n$ vectors in $\mathbb{Y}$. Hence, there exists $\lambda_{1},\dots,\lambda_{n}\in\mathbb{R}_+$ and $\mathbf{y}^{1},\dots,\mathbf{y}^{n}\in\mathbb{Y}$ such that $\mathbf{z}^*=\sum_{i=1}^n\lambda_i\mathbf{y}^i$. Therefore, problem $\mathscr{P}'$ is equivalent to:
$$\mathscr{P}':\min_{\substack{
\mathbf{z}=\sum_{i=1}^n\lambda_i\mathbf{y}^i\\
\lambda_{1},\dots,\lambda_{n}\in\mathbb{R}_+\\
\mathbf{y}^{1},\dots,\mathbf{y}^{n}\in\mathbb{Y}
}}  \|\mathbf{z}-\hat{\mathbf{z}}\|_2^2,$$
which is a significantly difficult problem because of the nonlinearity in the objective function and constraints, and the existence of $\mathcal{O}(n)$ continuous and $\mathcal{O}(n^2)$ discrete variables.

\section{Solution Methodology} \label{SectionSolutionMethod}

We apply the Frank-Wolfe Algorithm (FWA) \citep{FrankAndWolfe1956,DemyanovAndRubinov1970}, also known as the \emph{conditional gradient method}. In the literature, the FWA is used to solve an optimization problem with a convex function over a compact convex domain \citep{Jaggi2013}. We remark that problem $\mathscr{P}''$ has a non-convex objective function and non-compact feasible region; hence, the direct application of the FWA is not possible. However, we show how problem $\mathscr{P}''$ can be converted and become suitable for the application of the FWA.

Define $\mathcal{H}_{\hat{\mathbf{z}}}:=\{\mathbf{z}\in\mathbb{R}_+^n|\hat{\mathbf{z}}^T\mathbf{z}=\hat{\mathbf{z}}^T\hat{\mathbf{z}}\}$ and $\mathcal{FH}_{\hat{\mathbf{z}}}:=\{\mathbf{z}\in\mathcal{H}_{\hat{\mathbf{z}}}|\mathbf{z}\in\mathbb{Z}\}$. The following theorem shows that solving problem $\min_{\mathbf{z}\in\mathcal{FH}_{\hat{\mathbf{z}}}}\|\mathbf{z}-\hat{\mathbf{z}}\|_2^2$, one can obtain an optimal solution for $\mathscr{P}'$.

\begin{figure}[h]
\begin{center}
    \begin{subfigure}[b]{.49\textwidth}
        \centering
    \tdplotsetmaincoords{70}{120}
    \begin{tikzpicture}[tdplot_main_coords, scale = 1.1]

    \draw[opacity=0.2,thin,dashed](0,0,0)--(0,4,0)--(0,4,4)--(0,0,4)--(0,0,0)--cycle;
    \draw[opacity=0.2,thin,dashed](0,0,0)--(4,0,0)--(4,4,0)--(0,4,0)--(0,0,0)--cycle;
    \draw[opacity=0.2,thin,dashed](0,0,0)--(4,0,0)--(4,0,4)--(0,0,4)--(0,0,0)--cycle;
    \draw[opacity=0.2,thin,dashed](4,0,0)--(4,4,0)--(4,4,4)--(4,0,4)--(4,0,0)--cycle;
    \draw[opacity=0.2,thin,dashed](0,4,0)--(4,4,0)--(4,4,4)--(0,4,4)--(0,4,0)--cycle;
    \draw[opacity=0.2,thin,dashed](0,0,4)--(4,0,4)--(4,4,4)--(0,4,4)--(0,0,4)--cycle;

    \draw[fill=black] (0,4,0) circle (0pt) node[right]{$\mathbb{R}^3_+$};
        \draw[thick,->](0,0,0)--(3,0,0) node[below left]{};
        \draw[thick,->](0,0,0)--(0,3,0) node[below right]{};
        \draw[thick,->](0,0,0)--(0,0,3) node[left]{};

        \draw[thick,->,brown](0,0,0)--(1,1,0) node[below right]{$\hat{\mathbf{z}}$};

        \draw[fill=black] (0,2,0) circle (0pt) node[above right]{$\mathcal{H}_{\hat{\mathbf{z}}}$};
        \draw[fill=black] (0,2,3.5) circle (0pt) node[above right]{$\mathcal{FH}_{\hat{\mathbf{z}}}$};

        \filldraw[draw=black,fill=blue,opacity=0.15](2,0,0)--(0,2,0)--(0,2,3)--(1,1,2)--(2,0,4/3)--(2,0,0)--cycle;
        \filldraw[draw=black,fill=green,opacity=0.15](2,0,4)--(2,0,4/3)--(1,1,2)--(0,2,3)--(0,2,4)--(2,0,4)--cycle;
        \draw[fill=black] (1,1,2) circle (1pt) node[right]{$\mathbf{z}^1$};
        \draw[fill=black] (0,2,3) circle (1pt) node[right]{$\mathbf{z}^2$};
        \draw[fill=black] (4/3,2/3,2) circle (1pt) node[left]{$\mathbf{z}^3$};
        \draw[fill=black] (2,0,4/3) circle (1pt) node[left]{$\mathbf{z}^4$};
    \end{tikzpicture}
    \end{subfigure}
    \hfill
    \caption{Illustration of $\mathcal{H}_{\hat{\mathbf{z}}}$ and $\mathcal{FH}_{\hat{\mathbf{z}}}$. Note that $\mathcal{FH}_{\hat{\mathbf{z}}}$ is not bounded. \label{FH_zhatillustFig}}
\end{center}
\end{figure}
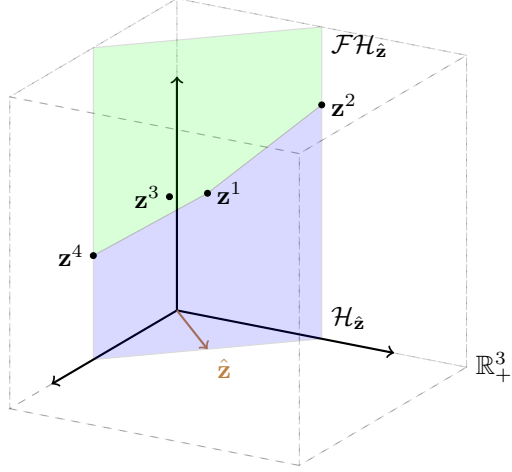

\begin{theorem}\label{pStar_pDStar_Theorem}
Let $\hat{\mathbf{z}}\in \mathbb{R}_+^n\backslash\{0\}$, $\mathbf{z}^{*}=\arg\min_{\mathbf{z}\in \mathbb{Z}}\|\mathbf{z}-\hat{\mathbf{z}}\|_2^2$, and $\mathbf{z}^{**}=\arg\min_{\mathbf{z}\in\mathcal{FH}_{\hat{\mathbf{z}}}}\|\mathbf{z}-\hat{\mathbf{z}}\|_2^2$. Then:

(a) $\mathbf{z}^{*}=0$ if and only if $\mathcal{FH}_{\hat{\mathbf{z}}}=\{\}$, and

(b) if $\mathcal{FH}_{\hat{\mathbf{z}}}\ne\{\}$, then
  $\mathbf{z}^{*}=\frac{\hat{\mathbf{z}}^T\hat{\mathbf{z}}}{\mathbf{z}^{**T}\mathbf{z}^{**}}\mathbf{z}^{**}$.
\end{theorem}
\noindent\textsc{Proof.}
(a, $\Rightarrow$) We first prove if $\mathbf{z}^{*}=0$, then $\mathcal{FH}_{\hat{\mathbf{z}}}=\{\}$, using a contradiction. Suppose that $\mathbf{z}^{*}=0$ and $\mathcal{FH}_{\hat{\mathbf{z}}}\ne\{\}$. Since $\mathcal{FH}_{\hat{\mathbf{z}}}$ is nonempty, then there exists $\bar{\mathbf{z}}\in\mathcal{FH}_{\hat{\mathbf{z}}}$. Hence, using the definition of $\mathcal{FH}_{\hat{\mathbf{z}}}$, we have $\bar{\mathbf{z}}\in\mathbb{R}^n_+$, $\hat{\mathbf{z}}^T\bar{\mathbf{z}}=\hat{\mathbf{z}}^T\hat{\mathbf{z}}$  and $\bar{\mathbf{z}}\in \mathbb{Z}$. Moreover, because $\hat{\mathbf{z}}^T\hat{\mathbf{z}}>0$, then $\bar{\mathbf{z}}\ne0$.

Define  $\bar{\bar{\mathbf{z}}}:=\frac{\hat{\mathbf{z}}^T\hat{\mathbf{z}}}{\bar{\mathbf{z}}^T\bar{\mathbf{z}}}\bar{\mathbf{z}}$. Since $\frac{\hat{\mathbf{z}}^T\bar{\mathbf{z}}}{\bar{\mathbf{z}}^T\bar{\mathbf{z}}}>0$ and $\bar{\mathbf{z}}\in \mathbb{Z}$, then $\bar{\bar{\mathbf{z}}}\in \mathbb{Z}$ (because $\mathbb{Z}$ is a cone). Since both $\hat{\mathbf{z}}$ and $\bar{\bar{\mathbf{z}}}$ are non-zero vectors in $\mathbb{R}^n_+$, then $\hat{\mathbf{z}}^T\bar{\bar{\mathbf{z}}}\ge0$. Moreover, $\hat{\mathbf{z}}$ and $\bar{\bar{\mathbf{z}}}$ are not orthogonal because $\bar{\bar{\mathbf{z}}}^T\hat{\mathbf{z}} =\frac{\hat{\mathbf{z}}^T\hat{\mathbf{z}}}{\bar{\mathbf{z}}^T\bar{\mathbf{z}}}\bar{\mathbf{z}}^T\hat{\mathbf{z}}>0$. In order to contradict the optimality of $\mathbf{z}^{*}=0$, we want to show the objective value of $\bar{\bar{\mathbf{z}}}$ is strictly better than that of $\mathbf{0}$. We have: $\|\hat{\mathbf{z}}-\bar{\bar{\mathbf{z}}}\|_2^2=\|\hat{\mathbf{z}}\|_2^2+\|\bar{\bar{\mathbf{z}}}\|_2^2-2\hat{\mathbf{z}}^T\bar{\bar{\mathbf{z}}}$.
Note that $\|\bar{\bar{\mathbf{z}}}\|_2^2-2\hat{\mathbf{z}}^T\bar{\bar{\mathbf{z}}}=
(\frac{\hat{\mathbf{z}}^T\hat{\mathbf{z}}}{\bar{\mathbf{z}}^T\bar{\mathbf{z}}})^2\bar{\mathbf{z}}^T\bar{\mathbf{z}}-2
(\frac{\hat{\mathbf{z}}^T\hat{\mathbf{z}}}{\bar{\mathbf{z}}^T\bar{\mathbf{z}}})\hat{\mathbf{z}}^T\bar{\mathbf{z}}= -(\frac{\hat{\mathbf{z}}^T\hat{\mathbf{z}}}{\bar{\mathbf{z}}^T\bar{\mathbf{z}}})\bar{\mathbf{z}}^T\bar{\mathbf{z}}<0$ (because $\hat{\mathbf{z}}^T\bar{\mathbf{z}}=\bar{\mathbf{z}}^T\bar{\mathbf{z}}$).  Therefore, $\|\hat{\mathbf{z}}-\bar{\bar{\mathbf{z}}}\|_2^2<\|\hat{\mathbf{z}}\|_2^2=\|\hat{\mathbf{z}}-\mathbf{0}\|_2^2$. This contradicts $\mathbf{z}^{*}=0$ because $\bar{\bar{\mathbf{z}}}$ is feasible and has a strictly better objective value. Hence, $\mathcal{FH}_{\hat{\mathbf{z}}}=\{\}$.

(a, $\Leftarrow$) We prove that if $\mathcal{FH}_{\hat{\mathbf{z}}}=\{\}$, then $\mathbf{z}^{*}=0$. Recall that in this paper we assume there exists $\mathbf{y}\in\mathbb{Y}$ such that $\mathbf{y}\ne0$. It follows that there exist $\mathbf{x}\in\mathbb{X}$ and $\mathbf{z}\in\mathbb{Z}$ such that $\mathbf{x}\ne0$ and $\mathbf{z}\ne0$. Since $\mathbb{Z}$ is a cone, then $\mathbf{0}\in\mathbb{Z}$. We must show that $\mathbf{0}$ is optimal ($\mathbf{0}$ has the smallest objective value among all $\mathbf{z}\in\mathbb{Z}$).

We claim that if $\bar{\mathbf{z}}\in \mathbb{Z}$, then $\hat{\mathbf{z}}^T\bar{\mathbf{z}}\le0$. We prove this claim using a contradiction. Let $\bar{\mathbf{z}}\in \mathbb{Z}$ such that $\hat{\mathbf{z}}^T\bar{\mathbf{z}}>0$. Define $\bar{\bar{\mathbf{z}}}:=\frac{\hat{\mathbf{z}}^T\hat{\mathbf{z}}}{\hat{\mathbf{z}}^T\bar{\mathbf{z}}}\bar{\mathbf{z}}$. Because $\frac{\hat{\mathbf{z}}^T\hat{\mathbf{z}}}{\hat{\mathbf{z}}^T\bar{\mathbf{z}}}>0$ and $\bar{\mathbf{z}}\in\mathbb{Z}$, then $\bar{\bar{\mathbf{z}}}\in \mathbb{Z}$. Moreover, $\hat{\mathbf{z}}^T\bar{\bar{\mathbf{z}}} =\hat{\mathbf{z}}^T(\frac{\hat{\mathbf{z}}^T\hat{\mathbf{z}}}{\hat{\mathbf{z}}^T\bar{\mathbf{z}}}\bar{\mathbf{z}})=\hat{\mathbf{z}}^T\hat{\mathbf{z}}$. Therefore, $\bar{\bar{\mathbf{z}}}\in\mathcal{FH}_{\hat{\mathbf{z}}}$, meaning that $\mathcal{FH}_{\hat{\mathbf{z}}}\ne\{\}$. This is a contradiction; hence, $\hat{\mathbf{z}}^T\bar{\mathbf{z}}\le0$, for all $\bar{\mathbf{z}}\in \mathbb{Z}$.

Let $\bar{\mathbf{z}}\in \mathbb{Z}$. Define $\dot{\mathbf{z}}:=(1-\frac{\hat{\mathbf{z}}^T\bar{\mathbf{z}}}{\hat{\mathbf{z}}^T\hat{\mathbf{z}}})\hat{\mathbf{z}}$ and $\ddot{\mathbf{z}}:=(\frac{\hat{\mathbf{z}}^T\bar{\mathbf{z}}}{\hat{\mathbf{z}}^T\hat{\mathbf{z}}}\hat{\mathbf{z}}-\bar{\mathbf{z}})$. Note that $\dot{\mathbf{z}}$ and $\ddot{\mathbf{z}}$ are orthogonal because
$\dot{\mathbf{z}}^T\ddot{\mathbf{z}}=(1-\frac{\hat{\mathbf{z}}^T\bar{\mathbf{z}}}{\hat{\mathbf{z}}^T\hat{\mathbf{z}}}) (\frac{\hat{\mathbf{z}}^T\bar{\mathbf{z}}}{\hat{\mathbf{z}}^T\hat{\mathbf{z}}}\hat{\mathbf{z}}^T\hat{\mathbf{z}}-\hat{\mathbf{z}}^T\bar{\mathbf{z}})=0$.
Moreover, we have $\|\dot{\mathbf{z}}\|_2^2\ge\|\hat{\mathbf{z}}\|_2^2$, because $(1-\frac{\hat{\mathbf{z}}^T\bar{\mathbf{z}}}{\hat{\mathbf{z}}^T\hat{\mathbf{z}}})\ge1$ (because $\hat{\mathbf{z}}^T\bar{\mathbf{z}}\le0$ and $\hat{\mathbf{z}}^T\hat{\mathbf{z}}>0$).
The objective value of $\bar{\mathbf{z}}$ is:
$\|\hat{\mathbf{z}}-\bar{\mathbf{z}}\|_2^2
=\|\hat{\mathbf{z}}-\frac{\hat{\mathbf{z}}^T\bar{\mathbf{z}}}{\hat{\mathbf{z}}^T\hat{\mathbf{z}}}\hat{\mathbf{z}} +\frac{\hat{\mathbf{z}}^T\bar{\mathbf{z}}}{\hat{\mathbf{z}}^T\hat{\mathbf{z}}}\hat{\mathbf{z}}-\bar{\mathbf{z}}\|_2^2
=\|\dot{\mathbf{z}}-\ddot{\mathbf{z}}\|_2^2
=\|\dot{\mathbf{z}}\|_2^2+\|\ddot{\mathbf{z}}\|_2^2
\ge\|\dot{\mathbf{z}}\|_2^2\ge\|\hat{\mathbf{z}}\|_2^2=\|\hat{\mathbf{z}}-\mathbf{0}\|_2^2$.
Therefore, we proved that $\|\hat{\mathbf{z}}-\bar{\mathbf{z}}\|_2^2\ge\|\hat{\mathbf{z}}-0\|_2^2$, for all $\bar{\mathbf{z}}\in \mathbb{Z}$. Thus, $\mathbf{z}^{*}=0$ is an optimal solution (there might be other optimal solutions).

(b) Since $\mathcal{FH}_{\hat{\mathbf{z}}}\ne\{\}$, using part (a), we have $\mathbf{z}^*\ne0$. We first claim that $\mathbf{z}^{*T}\hat{\mathbf{z}}>0$. Since $\mathbf{z}^{*},\ \hat{\mathbf{z}}\in\mathbb{R}^n_+$, then, $\mathbf{z}^{*T}\hat{\mathbf{z}}\ge0$. Suppose to the contrary that $\mathbf{z}^{*T}\hat{\mathbf{z}}=0$ (meaning that $\mathbf{z}^{*}$ and $\hat{\mathbf{z}}$ are orthogonal). Hence, $\|\hat{\mathbf{z}}-\mathbf{z}^{*}\|_2^2=\|\hat{\mathbf{z}}\|_2^2+\|\mathbf{z}^{*}\|_2^2>\|\hat{\mathbf{z}}\|_2^2=\|\hat{\mathbf{z}}-\mathbf{0}\|_2^2$. Since $\mathbf{0}\in \mathbb{Z}$ (see part (a)), this contradicts the optimality of $\mathbf{z}^{*}$, because $\mathbf{0}$ has a strictly better objective value.

We next show that $\mathbf{z}^*$ and $\mathbf{z}^{**}$ are unique. We prove this for $\mathbf{z}^*$ using a contradiction (the proof for $\mathbf{z}^{**}$ is similar). Suppose that $\dot{\mathbf{z}},\ddot{\mathbf{z}}\in\mathbb{Z}$ such that $\dot{\mathbf{z}}\ne\ddot{\mathbf{z}}$ and both $\dot{\mathbf{z}}$ and $\ddot{\mathbf{z}}$ are optimal for problem $\min_{\mathbf{z}\in \mathbb{Z}}\|\mathbf{z}-\hat{\mathbf{z}}\|_2^2$; hence $\|\dot{\mathbf{z}}-\hat{\mathbf{z}}\|_2^2=\|\ddot{\mathbf{z}}-\hat{\mathbf{z}}\|_2^2\le\|\mathbf{z}-\hat{\mathbf{z}}\|_2^2$, for all $\mathbf{z}\in\mathbb{Z}$. Define $\tilde{\mathbf{z}}:=\frac{1}{2}(\dot{\mathbf{z}}+\ddot{\mathbf{z}})$. Note that $\tilde{\mathbf{z}}\in\mathbb{Z}$ because $\mathbb{Z}$ is a convex set. We show that the objective value of $\tilde{\mathbf{z}}$ is strictly better than that of $\dot{\mathbf{z}}$ (or $\ddot{\mathbf{z}}$). We have:
$\|\tilde{\mathbf{z}}-\hat{\mathbf{z}}\|_2^2
=\|\frac{1}{2}(\dot{\mathbf{z}}+\ddot{\mathbf{z}})-\hat{\mathbf{z}}\|_2^2
=\frac{1}{4}\|\dot{\mathbf{z}}-\hat{\mathbf{z}}\|_2^2+\frac{1}{4}\|\ddot{\mathbf{z}}-\hat{\mathbf{z}}\|_2^2 +\frac{1}{2}(\dot{\mathbf{z}}-\hat{\mathbf{z}})^T(\ddot{\mathbf{z}}-\hat{\mathbf{z}})
\le\|\dot{\mathbf{z}}-\hat{\mathbf{z}}\|_2^2$. The equality holds only if the angle between $(\dot{\mathbf{z}}-\hat{\mathbf{z}})$ and $(\ddot{\mathbf{z}}-\hat{\mathbf{z}})$ is 0, in which case we must have $\dot{\mathbf{z}}=\ddot{\mathbf{z}}$ (because the lengthes of vectors $(\dot{\mathbf{z}}-\hat{\mathbf{z}})$ and $(\ddot{\mathbf{z}}-\hat{\mathbf{z}})$ are equivalent). This contradicts $\dot{\mathbf{z}}\ne\ddot{\mathbf{z}}$; hence $\mathbf{z}^*$ is unique. As we use the convexity of the feasible region to prove the uniqueness of $\mathbf{z}^*$, the proof for $\mathbf{z}^{**}$ is almost identical (because $\mathcal{FH}_{\hat{\mathbf{z}}}$ is a convex set).

Define $\lambda^*:=\frac{\hat{\mathbf{z}}^T\hat{\mathbf{z}}}{\mathbf{z}^{*T}\hat{\mathbf{z}}}$. Then, $\lambda^*>0$, because $\hat{\mathbf{z}}^T\hat{\mathbf{z}}>0$ and $\mathbf{z}^{*T}\hat{\mathbf{z}}>0$. Moreover, $\lambda^*\mathbf{z}^{*}\in\mathcal{FH}_{\hat{\mathbf{z}}}$, because $\lambda^*\mathbf{z}^{*}\in\mathbb{Z}$ and
$\hat{\mathbf{z}}^T(\lambda^*\mathbf{z}^{*}) =\frac{\hat{\mathbf{z}}^T\hat{\mathbf{z}}}{\mathbf{z}^{*T}\hat{\mathbf{z}}}\hat{\mathbf{z}}^T\mathbf{z}^{*}
=\hat{\mathbf{z}}^T\hat{\mathbf{z}}$. Moreover, $\frac{\mathbf{z}^{**}}{\lambda^*}\in\mathbb{Z}$, because $\frac{1}{\lambda^*}>0$, $\mathbf{z}^{**}\in\mathbb{Z}$, and $\mathbb{Z}$ is a cone. Hence, in summary, we have $\mathbf{z}^{**},\ \lambda^*\mathbf{z}^{*}\in\mathcal{FH}_{\hat{\mathbf{z}}}$ and $\frac{\mathbf{z}^{**}}{\lambda^*},\ \mathbf{z}^*\in\mathbb{Z}$.

Finally, we prove that $\mathbf{z}^{**}=\lambda^*\mathbf{z}^*$ using a contradiction. Suppose to the contrary that $\mathbf{z}^{**}\ne\lambda^*\mathbf{z}^{*}$ (or $\frac{1}{\lambda^*}\mathbf{z}^{**}\ne\mathbf{z}^{*}$). Since
$\mathbf{z}^*$ is unique and $\frac{1}{\lambda^*}\mathbf{z}^{**}\in\mathbb{Z}$, then we have:
\begin{eqnarray*}
  &&\|\mathbf{z}^{*}-\hat{\mathbf{z}}\|_2^2<\|\frac{1}{\lambda^*}\mathbf{z}^{**}-\hat{\mathbf{z}}\|^2_2 \\
  &&\Longrightarrow\|\mathbf{z}^{*}-\frac{1}{\lambda^*}\hat{\mathbf{z}}-(1-\frac{1}{\lambda^*})\hat{\mathbf{z}}\|_2^2 <\|\frac{1}{\lambda^*}\mathbf{z}^{**}-\frac{1}{\lambda^*}\hat{\mathbf{z}}-(1-\frac{1}{\lambda^*})\hat{\mathbf{z}}\|^2_2 \\
  &&\Longrightarrow\|\mathbf{z}^{*}-\frac{1}{\lambda^*}\hat{\mathbf{z}}\|_2^2+\|(1-\frac{1}{\lambda^*})\hat{\mathbf{z}}\|_2^2 <\|\frac{1}{\lambda^*}\mathbf{z}^{**}-\frac{1}{\lambda^*}\hat{\mathbf{z}}\|_2^2+\|(1-\frac{1}{\lambda^*})\hat{\mathbf{z}}\|^2_2 \\
  &&\Longrightarrow\|\mathbf{z}^{*}-\frac{1}{\lambda^*}\hat{\mathbf{z}}\|_2^2 <\|\frac{1}{\lambda^*}\mathbf{z}^{**}-\frac{1}{\lambda^*}\hat{\mathbf{z}}\|_2^2\\
  &&\Longrightarrow\|\lambda^*\mathbf{z}^{*}-\hat{\mathbf{z}}\|_2^2 <\|\mathbf{z}^{**}-\hat{\mathbf{z}}\|_2^2.
\end{eqnarray*}
The second line follows from the fact that $(\mathbf{z}^{*}-\frac{1}{\lambda^*}\hat{\mathbf{z}})$ and $(1-\frac{1}{\lambda^*})\hat{\mathbf{z}}$ are orthogonal, and $(\frac{1}{\lambda^*}\mathbf{z}^{**}-\frac{1}{\lambda^*}\hat{\mathbf{z}})$ and $(1-\frac{1}{\lambda^*})\hat{\mathbf{z}}$ are orthogonal. On the other hand, since
$\mathbf{z}^{**}$ is unique and $\lambda^*\mathbf{z}^{*}\in\mathcal{FH}_{\hat{\mathbf{z}}}$, then we must have: $\|\mathbf{z}^{**}-\hat{\mathbf{z}}\|_2^2<\|\lambda^*\mathbf{z}^*-\hat{\mathbf{z}}\|^2_2$. This is a contradiction; hence, $\mathbf{z}^{**}=\lambda^*\mathbf{z}^*$. Note that $\lambda^*=\frac{\hat{\mathbf{z}}^T\hat{\mathbf{z}}}{\mathbf{z}^{*T}\hat{\mathbf{z}}} =\frac{\mathbf{z}^{**T}\mathbf{z}^{**}}{\hat{\mathbf{z}}^T\hat{\mathbf{z}}}$. Thus, we showed that $\mathbf{z}^{*}=\frac{1}{\lambda^*}\mathbf{z}^{**}=\frac{\hat{\mathbf{z}}^T\hat{\mathbf{z}}}{\mathbf{z}^{**T}\mathbf{z}^{**}}\mathbf{z}^{**}$, and hence the proof is complete.
$\Box$

According to Theorem \ref{pStar_pDStar_Theorem}, we need to solve $\mathbf{z}^{**}=\arg\min_{\mathbf{z}\in\mathcal{FH}_{\hat{\mathbf{z}}}}\|\mathbf{z}-\hat{\mathbf{z}}\|_2^2$. We denote this problem by $\mathscr{P}'''$. If $\mathscr{P}'''$ is infeasible, then the optimal solution of $\mathscr{P}'$ is $\mathbf{z}^*=\mathbf{0}$; otherwise, $\mathbf{z}^{*}=\frac{\hat{\mathbf{z}}^T\hat{\mathbf{z}}}{\mathbf{z}^{**T}\mathbf{z}^{**}}\mathbf{z}^{**}$. In the remainder, we show how the FWA can be applied to solve $\mathscr{P}'''$. An initial feasible point is given $\mathbf{z}^{(0)}\in\mathcal{FH}_{\hat{\mathbf{z}}}$. At each iteration $k=0,1,2,\cdots$, we compute the gradient of the objective function at the current point $\mathbf{z}^{(k)}$, which is $\nabla f(\mathbf{z}^{(k)}):=(\mathbf{z}^{(k)}-\hat{\mathbf{z}})$, and then solve a minimization problem $\mathcal{M}(k): \min_{\mathbf{z}\in\mathcal{FH}_{\hat{\mathbf{z}}}}\mathbf{z}^T\nabla f(\mathbf{z}^{(k)})$. The current point for the next iteration, $\mathbf{z}^{(k+1)}$, is then updated. Problem $\mathcal{M}(k)$ is formulated as follows:
$$\mathcal{M}(k): \min_{\substack{\mathbf{y}\in\mathbb{Y},\ \lambda\ge0,\\
\hat{\mathbf{z}}^T(\lambda\mathbf{y})=\hat{\mathbf{z}}^T\hat{\mathbf{z}}}}
(\lambda\mathbf{y})^T\nabla f(\mathbf{z}^{(k)})$$

We remark that if $\mathbb{Y}$ is the set of feasible solutions of a BLP, then $\mathcal{M}(k)$ can be formulated as a mixed-integer linear programming (MILP) problem with $n$ binary and $n+1$ continuous variables. Moreover, if $\mathbb{Y}$ is given numerically, then the solution of  $\mathcal{M}(k)$ is obtained as $\min_{\mathbf{y}\in\mathbb{Y}:\hat{\mathbf{z}}^T\mathbf{y}\ne0}\{
(\frac{\hat{\mathbf{z}}^T\hat{\mathbf{z}}}{\hat{\mathbf{z}}^T\mathbf{y}})\mathbf{y}^T\nabla f(\mathbf{z}^{(k)})\}$, which can be performed in $\mathcal{O}(mn)$.

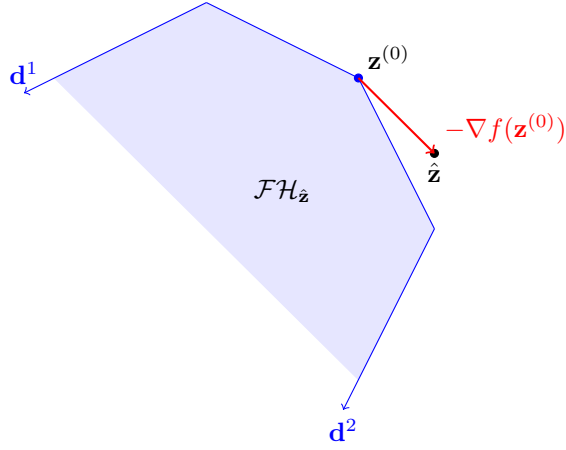
\begin{figure}[h]
\begin{center}
    \usetikzlibrary{arrows,intersections}
        \begin{tikzpicture}
          \coordinate (O) at (0,0);
          \filldraw[draw=blue,fill=blue,opacity=0.1] (-3,0) -- (-1,-1)-- (0,-3)-- (-1,-5)--(-5,-1)--(-3,0)--cycle;
          \draw[->,blue] (-3,0) -- (-5.4,-1.2) coordinate[label = {above:$\mathbf{d}^1$}] (d1);
          \draw[->,blue] (0,-3) -- (-1.2,-5.4) coordinate[label = {below:$\mathbf{d}^2$}] (d2);
          \draw[blue] (-3,0) -- (-1,-1)--(0,-3);
          \draw[draw=black,fill=black] (0,-2) circle (1.5pt) node[below]{ $\hat{\mathbf{z}}$};
          \draw[draw=blue,fill=blue] (-1,-1) circle (1.5pt) node[above right]{ $\mathbf{z}^{(0)}$};
          \draw[draw=blue,fill=blue] (-2.5,-2.5) circle (0pt) node[right]{ $\mathcal{FH}_{\hat{\mathbf{z}}}$};
          \draw[->,red,thick] (-1,-1)--(0,-2) coordinate[label = {above right:$-\nabla f(\mathbf{z}^{(0)})$}] (grad);
    \end{tikzpicture}
    \caption{An example to highlight the importance of the existence of $\mathcal{CFH}_{\hat{\mathbf{z}}}$.\label{compactFig}}
\end{center}
\end{figure}

The applicability of the FW depends on generating a feasible solution at each iteration. This is guaranteed if, for example, the feasible region of $\mathcal{M}(k)$ is compact. We remark that $\mathcal{FH}_{\hat{\mathbf{z}}}$ is a closed and convex polyhedron, but not necessarily compact because it can be unbounded (see Fig. \ref{FH_zhatillustFig}). We must ensure there exists a nonempty and compact polyhedron---denote by $\mathcal{CFH}_{\hat{\mathbf{z}}}$---that contains an optimal solution of $\mathcal{M}(k)$. This result is not obvious as illustrated in Fig. \ref{compactFig}. Suppose $\mathcal{FH}_{\hat{\mathbf{z}}}$, $\hat{\mathbf{z}}$, and $\mathbf{z}^{(0)}$ are as shown in Fig. \ref{compactFig}, and $\mathcal{FH}_{\hat{\mathbf{z}}}$ is unbounded with extreme directions $\mathbf{d}^1$ and $\mathbf{d}^2$. The objective value of $\mathcal{M}(0)$ improves in the direction of $-\nabla f(\mathbf{z}^{(0)})$. Hence, the objective value approaches to $-\infty$ in the direction of $\mathbf{d}^2$, and there exists no optimal solution. If this case happens, the FW algorithm cannot generate $\mathbf{z}^{(1)}$ and hence cannot proceed (A situation similar to Fig. \ref{compactFig} never happens as proven in Theorem \ref{compactTheorem}).

We first show that, if $\mathcal{FH}_{\hat{\mathbf{z}}}\ne\{\}$, then there exists a nonempty and compact polyhedron, which we denote by $\mathcal{CFH}_{\hat{\mathbf{z}}}$ ($\mathcal{C}$ stands for ``compact"), that contains an optimal solution of $\mathcal{M}(k)$. We also characterize the diameter of $\mathcal{CFH}_{\hat{\mathbf{z}}}$, which is denoted by $\mathcal{D}$. Define $\rho:=\max_{\mathbf{y}\in\mathbb{Y},\ \mathbf{y}^{T}\hat{\mathbf{z}}\ne0}
   \frac{\mathbf{y}^{T}\mathbf{y}}{(\mathbf{y}^{T}\hat{\mathbf{z}})^2}$.

\begin{theorem}\label{compactTheorem}
There exists a nonempty and compact polyhedron, denoted by $\mathcal{CFH}_{\hat{\mathbf{p}}}$, such that, for all $k\ge0$, the optimal solution of $\mathcal{M}(k)$ belongs to $\mathcal{CFH}_{\hat{\mathbf{p}}}$. Moreover, the diameter of $\mathcal{CFH}_{\hat{\mathbf{p}}}$ satisfies: $\mathcal{D}^2\le2\rho(\hat{\mathbf{z}}^T\hat{\mathbf{z}})^2$.
\end{theorem}

We remark that $\rho$ is bounded because $0\le\frac{\mathbf{y}^{T}\mathbf{y}}{(\mathbf{y}^{T}\hat{\mathbf{z}})^2}<+\infty$, for all $\mathbf{y}\in\mathbb{Y}$ such that $\mathbf{y}^{T}\hat{\mathbf{z}}\ne0$. Therefore, Theorem \ref{compactTheorem} provides a remedy for the unboundedness of the feasible region. The following Lemma characterizes the extreme points and extreme directions of $\mathcal{H}_{\hat{\mathbf{z}}}$, which will be used in the proof of Theorem \ref{compactTheorem} later in this section. Define $N=\{1,2,\dots,n\}$, $I_{\emptyset}:=\{i\in N|\hat{z}_i=0\}$, and $I_{1}:=\{i\in N|\hat{z}_i>0\}$.

\begin{lemma}\label{nonnegativePhatSimplex}
The vertices of $\mathcal{H}_{\hat{\mathbf{z}}}$ are $\frac{\hat{\mathbf{z}}^T\hat{\mathbf{z}}}{\hat{z}_i}\mathbf{e}^i$, for all $i\in I_{1}$, and the extreme directions of $\mathcal{H}_{\hat{\mathbf{z}}}$ are $\mathbf{e}^i$, for all $i\in I_{\emptyset}$, where $\mathbf{e}^i$ is the vector of all zeros except for $i$'th entry which is 1.
\end{lemma}
\noindent\textsc{Proof.} $\mathcal{H}_{\hat{\mathbf{z}}}$ is defined by 1 equality and $n$ inequalities; hence, the vertices of $\mathcal{H}_{\hat{\mathbf{z}}}$ are obtained by setting $n-1$ of inequalities to equalities and solving them together with the equation $\hat{\mathbf{z}}^T\mathbf{z}=\hat{\mathbf{z}}^T\hat{\mathbf{z}}$. Consider $i'$ such that $\hat{z}_{i'}>0$. Setting $\hat{z}_{i}=0$, for all $i\ne i'$, we obtain $\hat{z}_{i'}=\frac{\hat{\mathbf{z}}^T\hat{\mathbf{z}}}{\hat{z}_i}$; hence, an extreme point $\frac{\hat{\mathbf{z}}^T\hat{\mathbf{z}}}{\hat{z}_{i'}}\mathbf{e}^{i'}$ is found. Therefore, for each $i\in I_{1}$, there exists an extreme point $\frac{\hat{\mathbf{z}}^T\hat{\mathbf{z}}}{\hat{z}_{i}}\mathbf{e}^{i}$.

Extreme directions of $\mathcal{H}_{\hat{\mathbf{z}}}$ are the extreme points of the set $\{\mathbf{d}\in\mathbb{R}^n|\hat{\mathbf{z}}^T\mathbf{d}=0,\mathbf{d}\ge0,\mathbf{1}^T{\mathbf{d}}=1\}$, where $\mathbf{1}$ is the vector of appropriate size with all entries equal to 1. This results in the extreme directions $\mathbf{e}^i,\forall i\in I_{\emptyset}$.
$\Box$

\noindent\textsc{Proof of Theorem \ref{compactTheorem}.} Because $\mathcal{FH}_{\hat{\mathbf{p}}}$ is a nonempty polyhedron and $\mathcal{FH}_{\hat{\mathbf{p}}}\subseteq\mathbb{R}^n_+$, then $\mathcal{FH}_{\hat{\mathbf{p}}}$ does not contain a line, and has some (at least one) extreme points. Moreover, $\mathcal{FH}_{\hat{\mathbf{z}}}\subseteq\mathcal{H}_{\hat{\mathbf{z}}}$; hence, a direction for $\mathcal{FH}_{\hat{\mathbf{z}}}$ is also a direction for $\mathcal{H}_{\hat{\mathbf{z}}}$. Then, as a corollary of Lemma \ref{nonnegativePhatSimplex}, a direction of $\mathcal{FH}_{\hat{\mathbf{z}}}$ can be represented as a non-negative combination of the extreme directions of $\mathcal{H}_{\hat{\mathbf{z}}}$.

We first show that $\mathcal{M}(k)$ does not have unbounded optimal value. Suppose to the contrary that the optimal value of $\mathcal{M}(k)$ is unbounded. Let $\bar{\mathbf{z}}\in\mathcal{FH}_{\hat{\mathbf{z}}}$ be arbitrary. There must exist a direction $\mathbf{d}$ such that: $(\bar{\mathbf{z}}+\eta\mathbf{d})\in\mathcal{FH}_{\hat{\mathbf{z}}}$, for all $\eta\in\mathbb{R}_+$, and the objective value of $(\bar{\mathbf{z}}+\eta'\mathbf{d})$ is strictly less than that of $(\bar{\mathbf{z}}+\eta''\mathbf{d})$, if $\eta'>\eta''$. Vector $\mathbf{d}$ can be written as a non-negative combination of $\mathbf{e}^i$'s, $i\in I_{\emptyset}$; hence, there exists $\xi_i\ge0$, for all $i\in I_{\emptyset}$, such that: $\mathbf{d}=\sum_{i\in I_{\emptyset}}\xi_i \mathbf{e}^i$.
Let us compute the difference between the objective values of $(\bar{\mathbf{z}}+\eta'\mathbf{d})$ and $(\bar{\mathbf{z}}+\eta''\mathbf{d})$:
\begin{eqnarray*}
  (\mathbf{z}^{(k)}-\hat{\mathbf{z}})^T((\bar{\mathbf{z}}+\eta'\mathbf{d})-(\bar{\mathbf{z}}+\eta''\mathbf{d}))&=&
  (\eta'-\eta'')(\mathbf{z}^{(k)}-\hat{\mathbf{z}})^T\mathbf{d} \\
  &=&   (\eta'-\eta'')(\mathbf{z}^{(k)}-\hat{\mathbf{z}})^T\sum_{i\in I_{\emptyset}}\xi_i \mathbf{e}^i \\
  &=&(\eta'-\eta'')\sum_{i\in I_{\emptyset}}\xi_iz^{(k)}_i\ge0,
\end{eqnarray*}
where, the last line follows from $\hat{z}_i=0$, for all $i\in I_{\emptyset}$ (also note that $\eta'>\eta''$). This is a contradiction because the objective value of $(\bar{\mathbf{z}}+\eta'\mathbf{d})$ is \textbf{not} strictly less than that of $(\bar{\mathbf{z}}+\eta''\mathbf{d})$.
Hence, $\mathcal{M}(k)$ has a bounded optimal value.

Additionally, since $\mathcal{FH}_{\hat{\mathbf{p}}}$ has at least one extreme point and does not contain lines, an extreme point of $\mathcal{FH}_{\hat{\mathbf{p}}}$ must be optimal. This extreme point is in the form of $\lambda\mathbf{y}$, where $0<\lambda<+\infty$ and $\mathbf{y}\in\mathbb{Y}$. We remark that there might be other alternative optimal solutions in the same form. However, there cannot be an optimal solution in the form of $\lambda\mathbf{y}$, where $\lambda=+\infty$ and $\mathbf{y}\in\mathbb{Y}$. Hence, let $\mathcal{CFH}_{\hat{\mathbf{p}}}$ denote the convex hall of all $\lambda\mathbf{y}$'s, where $0<\lambda<+\infty$ and $\mathbf{y}\in\mathbb{Y}$. Thus, we proved that the optimal solution of $\mathcal{M}(k)$ belongs to $\mathcal{CFH}_{\hat{\mathbf{z}}}$. Moreover, $\mathcal{CFH}_{\hat{\mathbf{z}}}$ is a nonempty and compact polyhedron.

We next obtain an upper bound on the diameter of $\mathcal{CFH}_{\hat{\mathbf{z}}}$. We have:
\begin{eqnarray*}
  \mathcal{D}^2 &=& \max_{\mathbf{z}^1,\mathbf{z}^2\in \mathcal{CFH}_{\hat{\mathbf{z}}}}\|\mathbf{z}^1-\mathbf{z}^2\|_2^2 \\
   &=& \max_{\substack{\lambda_i\mathbf{y}^{iT}\hat{\mathbf{z}}=\hat{\mathbf{z}}^T\hat{\mathbf{z}},\ i=1,2,\\ \lambda_1,\lambda_2>0,\ \mathbf{y}^1,\mathbf{y}^2\in\mathbb{Y}}}\|\lambda_1\mathbf{y}^1-\lambda_2\mathbf{y}^2\|_2^2\\
   &=& \max_{\substack{\lambda_i\mathbf{y}^{iT}\hat{\mathbf{z}}=\hat{\mathbf{z}}^T\hat{\mathbf{z}},\ i=1,2,\\ \lambda_1,\lambda_2>0,\ \mathbf{y}^1,\mathbf{y}^2\in\mathbb{Y}}} \lambda_1^2\mathbf{y}^{1T}\mathbf{y}^1+\lambda_2^2\mathbf{y}^{2T}\mathbf{y}^2-2\lambda_1\lambda_2\mathbf{y}^{1T}\mathbf{y}^2\\
   &\le& 2\max_{\substack{\lambda\mathbf{y}^{T}\hat{\mathbf{z}}=\hat{\mathbf{z}}^T\hat{\mathbf{z}},\\ \lambda>0,\ \mathbf{y}\in\mathbb{Y}}} \lambda^2\mathbf{y}^{T}\mathbf{y}\\
   &=& 2\max_{\mathbf{y}\in\mathbb{Y},\ \mathbf{y}^{T}\hat{\mathbf{z}}\ne0}
   \big(\frac{\hat{\mathbf{z}}^T\hat{\mathbf{z}}}{\mathbf{y}^{T}\hat{\mathbf{z}}}\big)^2\mathbf{y}^{T}\mathbf{y}\\
   &=& 2(\hat{\mathbf{z}}^T\hat{\mathbf{z}})^2\max_{\mathbf{y}\in\mathbb{Y},\ \mathbf{y}^{T}\hat{\mathbf{z}}\ne0}
   \frac{\mathbf{y}^{T}\mathbf{y}}{(\mathbf{y}^{T}\hat{\mathbf{z}})^2}
\end{eqnarray*}
The first line follows from the definition of $\mathcal{D}$. In the second line, we use the fact that $\mathbf{z}^1$ and $\mathbf{z}^2$ must be extreme points of $\mathcal{CFH}_{\hat{\mathbf{z}}}$; hence, they should be in the form of $\mathbf{z}^i=\lambda_i\mathbf{y}^i$ where $\lambda^i>0$, $\mathbf{y}^i\in \mathbb{Y}$, and $\lambda_i\mathbf{y}^i$ must be on the hyperplane $\mathbf{z}^T\hat{\mathbf{z}}=\hat{\mathbf{z}}^T\hat{\mathbf{z}}$. In the forth line, we eliminate $-2\lambda_1\lambda_2\mathbf{y}^{1T}\mathbf{y}^2$ because it is always non-positive. As a result, the problem decomposes into two identical problems. In the fifth line, we substitute $\lambda$ with $\frac{\hat{\mathbf{z}}^T\hat{\mathbf{z}}}{\mathbf{y}^{T}\hat{\mathbf{z}}}$ but we have to ensure that the denominator is non-zero. Hence, the proof is complete.
$\Box$

Due to Theorem \ref{compactTheorem}, our problem is equivalent to $\min_{\mathbf{z}\in\mathcal{CFH}_{\hat{\mathbf{z}}}}\|\mathbf{z}-\hat{\mathbf{z}}\|_2^2$, where the objective function is convex and the feasible region is nonempty and compact. It is well-known (see for example \citet{Jaggi2013}) that for each $k\ge1$, the iterates $\mathbf{z}^{(k)}$ of the FWA satisfy:
$$\|\mathbf{z}^{(k)}-\hat{\mathbf{z}}\|^2_2-\|\mathbf{z}^{**}-\hat{\mathbf{z}}\|^2_2\le\frac{2\mathcal{D}^2}{k+2},$$
where $\mathbf{z}^{**}$ is the optimal solution of $\min_{\mathbf{z}\in\mathcal{CFH}_{\hat{\mathbf{z}}}}\|\mathbf{z}-\hat{\mathbf{z}}\|_2^2$.
Combining this result with Theorem \ref{compactTheorem}, we obtain:
$$\|\mathbf{z}^{(k)}-\hat{\mathbf{z}}\|^2_2-\|\mathbf{z}^{**}-\hat{\mathbf{z}}\|^2_2\le\frac{4\rho(\hat{\mathbf{z}}^T\hat{\mathbf{z}})^2}{k+2}.$$

Hence, the proposed Frank-Wolfe based approach solves our problem at a linear convergence rate, i.e. the optimization error after $k$ iterations will decrease with $\mathcal{O}(\frac{1}{k})$.

\subsection{Example}

We graphically show the steps of applying our approach to the example presented in Figs. \ref{RCPillustFig} and \ref{FH_zhatillustFig}. The problem is solved in two FWA iterations. In iteration 0, we start from $\mathbf{z}^{(0)}$, find $\nabla f(\mathbf{z}^{(0)})$, and maximize in the direction of $-\nabla f(\mathbf{z}^{(0)})$ to obtain $\tilde{\mathbf{z}}$. Then, $\mathbf{z}^{(1)}$ is found as the closest point to $\hat{\mathbf{z}}$ in the convex combination of $\mathbf{z}^{(0)}$ and $\tilde{\mathbf{z}}$. In iteration 1, we find $\nabla f(\mathbf{z}^{(0)})$, and then obtain $\tilde{\mathbf{z}}$, and the algorithm stops.

\begin{figure}[h]
\begin{center}
    \begin{subfigure}[b]{.49\textwidth}
        \centering
    \tdplotsetmaincoords{70}{120}
    \begin{tikzpicture}[tdplot_main_coords, scale = 1.1]

    \draw[opacity=0.2,thin,dashed](0,0,0)--(0,4,0)--(0,4,4)--(0,0,4)--(0,0,0)--cycle;
    \draw[opacity=0.2,thin,dashed](0,0,0)--(4,0,0)--(4,4,0)--(0,4,0)--(0,0,0)--cycle;
    \draw[opacity=0.2,thin,dashed](0,0,0)--(4,0,0)--(4,0,4)--(0,0,4)--(0,0,0)--cycle;
    \draw[opacity=0.2,thin,dashed](4,0,0)--(4,4,0)--(4,4,4)--(4,0,4)--(4,0,0)--cycle;
    \draw[opacity=0.2,thin,dashed](0,4,0)--(4,4,0)--(4,4,4)--(0,4,4)--(0,4,0)--cycle;
    \draw[opacity=0.2,thin,dashed](0,0,4)--(4,0,4)--(4,4,4)--(0,4,4)--(0,0,4)--cycle;

    \draw[fill=black] (0,4,0) circle (0pt) node[right]{$\mathbb{R}^3_+$};
        \draw[thick,->](0,0,0)--(3,0,0) node[below left]{};
        \draw[thick,->](0,0,0)--(0,3,0) node[below right]{};
        \draw[thick,->](0,0,0)--(0,0,3) node[left]{};

        \draw[thick,->,brown](0,0,0)--(1,1,0) node[below right]{$\hat{\mathbf{z}}$};

        \draw[fill=black] (0,2,3.5) circle (0pt) node[above right]{$\mathcal{FH}_{\hat{\mathbf{z}}}$};

        \filldraw[draw=black,fill=green,opacity=0.15](2,0,4)--(2,0,4/3)--(1,1,2)--(0,2,3)--(0,2,4)--(2,0,4)--cycle;
        \draw[fill=black] (1,1,2) circle (1pt) node[right]{$\mathbf{z}^{(0)}$};
        \draw[fill=black] (2,0,4/3) circle (1pt) node[left]{$\tilde{\mathbf{z}}$};
        \draw[thick,->,red](1,1,2)--(1,1,0) node[above right, yshift=.8cm]{$-\nabla f(\mathbf{z}^{(0)})$};
        \draw[fill=black] (17/11,5/11,54/33) circle (1pt) node[above]{$\mathbf{z}^{(1)}$};
        \draw[fill=black] (5,4,0) circle (0pt) node[below]{Iteration 0: $\tilde{\mathbf{z}}$ minimizes $\mathcal{M}(0)$ };
    \end{tikzpicture}
    \end{subfigure}
        \begin{subfigure}[b]{.49\textwidth}
        \centering
    \tdplotsetmaincoords{70}{120}
    \begin{tikzpicture}[tdplot_main_coords, scale = 1.1]

    \draw[opacity=0.2,thin,dashed](0,0,0)--(0,4,0)--(0,4,4)--(0,0,4)--(0,0,0)--cycle;
    \draw[opacity=0.2,thin,dashed](0,0,0)--(4,0,0)--(4,4,0)--(0,4,0)--(0,0,0)--cycle;
    \draw[opacity=0.2,thin,dashed](0,0,0)--(4,0,0)--(4,0,4)--(0,0,4)--(0,0,0)--cycle;
    \draw[opacity=0.2,thin,dashed](4,0,0)--(4,4,0)--(4,4,4)--(4,0,4)--(4,0,0)--cycle;
    \draw[opacity=0.2,thin,dashed](0,4,0)--(4,4,0)--(4,4,4)--(0,4,4)--(0,4,0)--cycle;
    \draw[opacity=0.2,thin,dashed](0,0,4)--(4,0,4)--(4,4,4)--(0,4,4)--(0,0,4)--cycle;

    \draw[fill=black] (0,4,0) circle (0pt) node[right]{$\mathbb{R}^3_+$};
        \draw[thick,->](0,0,0)--(3,0,0) node[below left]{};
        \draw[thick,->](0,0,0)--(0,3,0) node[below right]{};
        \draw[thick,->](0,0,0)--(0,0,3) node[left]{};

        \draw[thick,->,brown](0,0,0)--(1,1,0) node[below right]{$\hat{\mathbf{z}}$};

        \draw[fill=black] (0,2,3.5) circle (0pt) node[above right]{$\mathcal{FH}_{\hat{\mathbf{z}}}$};

        \filldraw[draw=black,fill=green,opacity=0.15](2,0,4)--(2,0,4/3)--(1,1,2)--(0,2,3)--(0,2,4)--(2,0,4)--cycle;
        \draw[fill=black] (1,1,2) circle (1pt) node[right]{$\mathbf{z}^{(0)}$};
        \draw[fill=black] (2,0,4/3) circle (1pt) node[left]{$\tilde{\mathbf{z}}$};
        \draw[thick,->,red](17/11,5/11,54/33)--(1,1,0) node[above, yshift=.8cm, xshift=.4cm]{$-\nabla f(\mathbf{z}^{(1)})$};
        \draw[fill=black] (17/11,5/11,54/33) circle (1pt) node[above]{$\mathbf{z}^{(1)}$};
        \draw[fill=black] (5,4,0) circle (0pt) node[below]{Iteration 1: $\tilde{\mathbf{z}}$ minimizes $\mathcal{M}(1)$ };
    \end{tikzpicture}
    \end{subfigure}
    \hfill
    \caption{Graphical illustration of the first two iterations of the FWA. \label{FWA_illustFig}}
\end{center}
\end{figure}
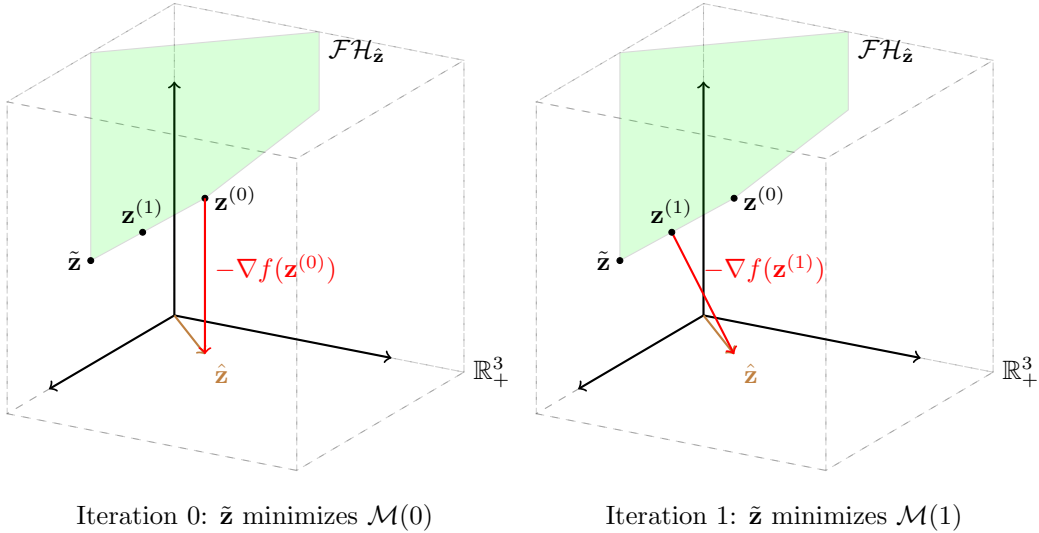

\section{Numerical Results}

We test the effectiveness of the FWA on a limited set of problem instances generated randomly. We take $|\mathbb{Y}|=1000$, and $n$ varies between 100 and 1000 in the increments of 100. We consider 1, 2, and 5 clusters for each problem size. Let NC denote the number of clusters and CF denote cluster coefficient---a measure of the closeness of the $\mathbf{y}$'s in a cluster. If CF=0, then all $\mathbf{y}$'s are randomly distributed in $\mathbb{R}^n_+$. In this case, we create each $\mathbf{y}$ by randomly generating $n$ nonnegative numbers. If CF$>$0, then we assume that the data points from $(\frac{(\ell-1)|\mathbb{Y}|}{NC}+1)$ to $(\frac{\ell|\mathbb{Y}|}{NC})$ are in the $\ell$'th cluster and the $(\frac{(\ell-1)|\mathbb{Y}|}{NC}+1)$'th data point is the center of the $\ell$'th cluster, for all $\ell=1,\dots,NC$. The center of cluster $\ell$, denoted by $\mathbf{y}^{c\ell}$, is created randomly. Then, a new point is generated in cluster $\ell$ using the formula $(CF\mathbf{y}^{c\ell}+\dot{\mathbf{y}})$, where $\dot{\mathbf{y}}$ is generated randomly in $\mathbb{R}^n_+$. Obviously, if CF=0, then there is no clustering, and as CF increases the angle between the points and the center of the cluster becomes smaller. Note that the lengths of $\mathbf{y}$'s are not important as we are interested in the rays generated by $\mathbf{y}$'s.

\begin{figure}
  \centering
  \includegraphics[width=10cm]{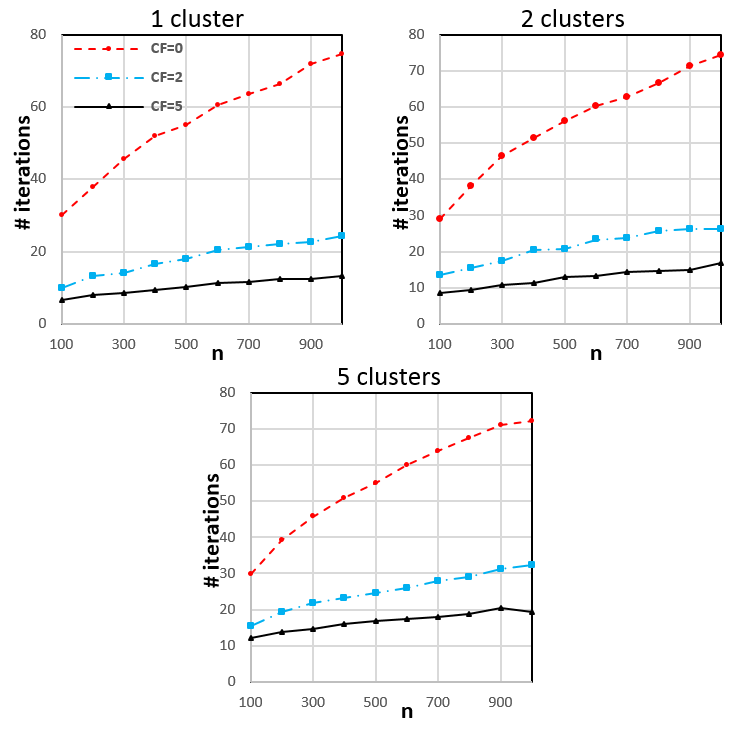}
  \caption{Effect of $n$, number of clusters, and cluster coefficient on the number of iterations.}\label{FW_CompFig}
\end{figure}

Fig. \ref{FW_CompFig} summarizes the result of applying our approach for 1, 2, and 5 clusters, where the horizontal and vertical axis show the number of dimensions and the number of iterations. Each point indicates the average of 50 randomly generated instances. Fig. \ref{FW_CompFig} shows that, as expected, the number of iterations (almost linearly) increases in $n$. As CF increases, the number of iterations decreases. Finally, the number of iterations increases in the number of clusters.

\section{Conclusions}

In this paper, we described a non-convex optimization problem with a noncompact-convex domain. This class of problems are found in manufacturing systems, clustering, machine learning, and statistics. We show how the FWA can be applied to these new class of problems by first proposing an equivalent problem with a convex objective function over a convex and non-compact domain, and then finding a compact set that contains the iterates of the FWA. A numerical example illustrates the steps of our approach. Finally, the numerical experiment shows the impact of $n$, number of clusters, and cluster coefficient on the number of iterations.

\end{document}